\documentclass{article}
\usepackage{amssymb}
\usepackage{amsmath}
\usepackage{amsfonts}
\usepackage{graphicx}

\newtheorem{theorem}{Theorem}[section]
\newtheorem{lemma}[theorem]{Lemma}

\newtheorem{remark}[theorem]{Remark}

\newcommand{\qed}{\hfill \rule{2.3mm}{2.3mm}}

\newcommand{\al}{\alpha}

\newcommand{\ve}{\varepsilon}
\newcommand{\vp}{\varphi}

\newcommand{\bee}{\begin{equation}}
\newcommand{\ee}{\end{equation}}

\newcommand{\R}{\mathbb{R}}
\newcommand{\N}{\mathbb{N}}

\newcommand{\Z}{\mathbb{Z}}
\newcommand{\BB}{\mathbb{B}}

\newcommand{\x}{x/\varepsilon}
\newcommand{\y}{y/\varepsilon}

\newcommand{\W}{{\cal W}}
\newcommand{\V}{{\cal V}}
\newcommand{\F}{{\cal F}}
\newcommand{\U}{{\cal U}}
\newcommand{\G}{{\cal G}}

\newcommand{\reff}[1]{(\ref{#1})}   
\begin{document}

\title{\bf
A Common Approach to Singular Perturbation and Homogenization I: Quasilinear ODE Systems}
%\footnote{\bf This paper appeared as preprint arXiv:2309.15611 in September 2024.
%Unfortunately Nikolai N. Nefedov passed away in October 2024.}
% is dedicated to the memory of Nikolai N. Nefedov, the coauthor of parts I and II of the present %paper, a great mathematician and a wonderful friend.
%He passed away in October 2024.}
%}

\newcounter{thesame}
\setcounter{thesame}{1}
\author{
Nikolai N. Nefedov \thanks{Lomonosov Moscow State University, Faculty of Physics, Vorob'jovy Gory Moscow, 10899, Russia. Nikolai N. Nefedov passed away in October 2024.
}
\ \ \ \ Lutz Recke \thanks{Humboldt University of Berlin, Institute of Mathematics, Rudower Chaussee 25, 12489 Berlin, Germany.
		{\small   E-mail:
			{\tt lutz.recke@hu-berlin.de}}
}}

%\setcounter{thesame}{2}
%\author{
%N%ikolai N. Nefedov \thanks{Lomonosov Moscow State University, Faculty of Physics, Vorob'jovy Gory, Moscow, 10899, Russia.
%}}

\date{}

\maketitle

\begin{abstract}
\noindent
We consider  periodic homogenization of
boundary value problems for quasilinear 
second-order ODE systems in divergence form of the type
$$
a(x,\x,u(x),u'(x))'= f(x,\x,u(x),u'(x)) \mbox{ for } x \in [0,1].
$$
For small $\varepsilon>0$ we show existence of solutions $u=u_\varepsilon$ as well as their local uniqueness close to 
%for $\|u-u_0\|_\infty \approx 0$,
%where $u_0$ is 
a given
non-degenerate solution $u_0$ to the homogenized problem, and we describe the rate of convergence to zero for $\varepsilon \to 0$ of the homogenization error $\|u_\varepsilon-u_0\|_\infty$. In particular, we show that this rate depends on the smoothness of the maps $a(\cdot,y,u,u')$ and $f(\cdot,y,u,u')$. 

Our assumptions are, roughly speaking, as follows: The maps $a,f:[0,1]\times\mathbb{R}\times\mathbb{R}^n\times\mathbb{R}^n\to\mathbb{R}^n$ are continuous, the maps $a(x,y,\cdot,\cdot)$ and $f(x,y,\cdot,\cdot)$ are $C^1$-smooth,
the maps $a(x,\cdot,u,u')$ and $f(x,\cdot,u,u')$ are 1-periodic, and the maps $a(x,y,u,\cdot)$ are strongly monotone and Lipschitz continuous  
uniformly with respect to $x$, $y$ and bounded $u$. No
global solution uniqueness is supposed.
Because $x$ is one-dimensional, neither correctors nor cell problems are needed. 

The main tool of the proofs is an abstract result of implicit function theorem type
which in the past has been applied to singularly perturbed nonlinear ODEs and elliptic and parabolic PDEs and, hence, which permits a common approach to
existence and local uniqueness results for singularly perturbed problems and and for homogenization problems.
\end{abstract}

\section{Introduction}
\setcounter{equation}{0}
\setcounter{theorem}{0}

In this paper we present an abstract result of implicit function theorem type (see Section \ref{secabstract}),
which in the past has been applied
in \cite{Butetc,But2022,Fiedler,NURS,
OmelchenkoRecke2015,Recke2022,
ReckeOmelchenko2008}  to singularly perturbed nonlinear ODEs and PDEs, in Parts II \cite{II}
and III \cite{III} to periodic homogenization of semilinear elliptic PDE systems and in \cite{AA2025} to non-periodic homogenization of semilinear ODE systems.
In the present paper we apply it to describe periodic homogenization for systems of quasilinear second-order ODEs in divergence form of the type
\bee
\label{ODE}
a(x,\x,u(x),u'(x))'= f(x,\x,u(x),u'(x))
\mbox{ for } x \in [0,1]
\ee
with one Dirichlet and one natural boundary condition
\bee
\label{bc}
u(0)=a(1,1/\ve,u(1),u'(1))=0
\ee
as well as with  other boundary conditions (see Section \ref{sec:other}).
Here $\ve>0$ is the small homogenization parameter, and we look for vector valued solutions $u:[0,1] \to \R^n$. The coefficient functions $a,f:[0,1]\times \R \times \R^n\times \R^n \to \R^n$
are supposed to be 1-periodic with respect to the second argument, i.e.
\bee
\label{perass}
a(x,y+1,u,u')=a(x,y,u,u')\mbox{ and }  
f(x,y+1,u,u')=f(x,y,u,u')
\ee
for all 
$x\in [0,1]$, $y\in \R$ and $u,u'\in\R^n$.
Further, we suppose that the maps $a$ and $f$ are continuous and that their first partial derivatives with respect to 
the third and fourth arguments exist and
are continuous, i.e.
\bee
\label{diffass}
a,\partial_ua,\partial_{u'}a,
f,\partial_uf,\partial_{u'}f
\mbox{ are continuous on }
[0,1]\times\R\times\R^n\times\R^n.
\ee
And finally, we suppose that the maps $a(x,y,u,\cdot)$ are strongly monotone and Lipschitz continuous 
uniformly with respect to $x$, $y$ and bounded $u$, i.e. that there exist constants $M \ge m >0$ such that for all 
$x \in [0,1]$, $y \in \R$
and $u,u_1',u_2' \in \R^n$
with $\|u\|\le 1$ we have
\bee
\label{mon}
\left.
\begin{array}{l}
\big(a(x,y,u,u_1')-
a(x,y,u,u_2')\big)\cdot(u_1'-u_2') \ge 
m\|u_1'-u_2'\|^2,\\
\|a(x,y,u,u_1')-a(x,y,u,u_2')\|
\|\le M\|u'_1-u'_2\|.
\end{array}
\right\}
\ee
Here and in what follows we denote by $v\cdot w$ the Euclidean scalar product  of vectors $v,w \in \R^n$, and $\|v\|:=\sqrt{v\cdot v}$ is the Euclidean norm of the vector $v \in \R^n$.

From assumption \reff{mon} it follows that for all $x\in [0,1]$, $y \in \R$ and $u \in \R^n$ with $\|u\|\le 1$ the maps $a(x,y,u,\cdot)$ are bijective from $\R^n$ onto $\R^n$. We denote 
\bee
\label{bdef}
b(x,y,u,\cdot):=
a(x,y,u,\cdot)^{-1}
%, \mbox{ i.e. }
%a(x,y,u,b(x,y,u,u'))=b(x,y,u,a(x,y,u,u'))=u'
\mbox{ and }
%\label{bnulldef}
b_0(x,u,u'):=\int_0^1b(x,y,u,u')dy
\ee 
for all $x \in [0,1]$, $y \in \R$
and $u,u'\in \R^n$
with $\|u\|\le 1$.
Also the maps $b_0(x,u,\cdot)$ are strongly monotone and Lipschitz continuous and, hence, bijective from $\R^n$ onto $\R^n$, and we denote 
\bee
\label{anulldef}
a_0(x,u,\cdot):=
b_0(x,u,\cdot)^{-1}
\ee
and
\bee
\label{fnulldef}
f_0(x,u,u'):=
\int_0^1f(x,y,u,b(x,y,u,a_0(x,u,u')))dy
\ee
for $x \in [0,1]$
and $u,u'\in\R^n$ with $\|u\|\le 1$,
and the boundary value problem
\bee
\label{hombvp}
a_0(x,u(x),u'(x))'=f_0(x,u(x),u'(x)) \mbox{ for } x \in [0,1],\;
u(0)=a_0(1,u(1),u'(1))=0
\ee
is, by definition, the homogenized version of the boundary value problem \reff{ODE}-\reff{bc}.

A vector function $u\in C^1([0,1];\R^n)$ is called weak solution to \reff{ODE}-\reff{bc} if it
satisfies the Dirichlet boundary condition $u(0)=0$ and
the variational equation
\bee
\label{vareq}
\left.
\begin{array}{r}
\displaystyle\int_0^1\Big(a(x,\x,u(x),u'(x))\cdot\vp'(x)
+f(x,\x,u(x),u'(x))\cdot\vp(x)\Big)dx=0\\
\mbox{ for all } \vp \in C^1([0,1];\R^n)
\mbox{ with } \vp(0)=0,
\end{array}
\right\}
\ee
and similar for the homogenized boundary value problem \reff{hombvp}
and its linearization \reff{linbvp}.
Weak solutions to \reff{ODE}-\reff{bc} are not classical solutions, i.e. they are not $C^2$-smooth, in general, because the maps $a(\cdot,\cdot,u,u')$ 
are not $C^1$-smooth, in general.

Now we formulate our result about existence and local uniqueness of
weak solutions $u=u_\ve$ to \reff{ODE}-\reff{bc} with $\ve \approx 0$, which are close to a given non-degenerate weak solution $u=u_0$ to \reff{hombvp}, and about the rate of convergence to zero for $\ve \to 0$ of the homogenization error $\|u_\ve-u_0\|_\infty$.
Here and in what follows we denote by
$\|u\|_\infty:=\max\{\|u(x)\|:\;
x \in [0,1]\}$
the maximum norm in the function space
$C([0,1];\R^n)$. 
\begin{theorem} 
\label{main}
Suppose \reff{perass}-\reff{mon}, and
let $u=u_0$ be a weak solution to the homogenized boundary value problem \reff{hombvp} such that $\|u_0\|_\infty<1$
and that the linearized homogenized boundary value problem 
\bee
\label{linbvp}
\left.
\begin{array}{l}
\Big(\partial_ua_0(x,u_0(x),u_0'(x))u(x)+\partial_{u'}a_0(x,u_0(x),u_0'(x))
u'(x)\Big)'\\
= \partial_uf_0(x,u_0(x),u_0'(x))u(x)
+\partial_{u'}f_0(x,u_0(x),u_0'(x))u'(x)
\mbox{ for }
x \in [0,1],\\
u(0)=\partial_ua_0(1,u_0(1),u_0'(1))u(1)+\partial_{u'}a_0(1,u_0(1),u_0'(1))u'(1)=0
\end{array}
\right\}
\ee
does not have weak solutions $u\not=0$.
Then the following is true:

(i) There exist $\ve_0>0$ and  $\delta>0$
such that for all $\ve \in (0,\ve_0]$ there exists exactly one 
weak solution $u=u_\ve$ 
to \reff{ODE}-\reff{bc} with $\|u(x)-u_0(x)\|+\|a(x,\x,u(x),u'(x))-a_0(x,u_0(x),u'_0(x))\| \le \delta$
for all $x \in [0,1]$. Moreover,
\bee
\label{est}
\|u_\ve(x)-u_0(x)\|
+\|a(x,\x,u_\ve(x),u_\ve'(x))-a_0(x,u_0(x),u'_0(x))\|
\to 0 \mbox{ for } \ve \to 0
\ee
uniformly with respect to $x \in [0,1]$.

(ii) For $\ve\in (0,\ve_0], x \in [0,1]$ denote $v_\ve(x):=a(x,\x,u_\ve(x),u_\ve'(x)), v_0(x):=
a_0(x,u_0(x),u'_0(x))$. Further, suppose that
\bee
\label{diffass1}
\partial_xa,
\partial_xf
\mbox{ exist and are continuous on }
[0,1]\times\R\times\R^n\times \R^n.
\ee
Then
\bee
\label{est1}
\|u_\ve-u_0\|_\infty+\|v_\ve-v_0\|_\infty
=O(\ve) \mbox{ for } \ve \to 0.
\ee
\end{theorem}

\begin{remark}
If $f(x.y,u,u')$ does not depend on $u'$, then for any weak solution $u$ to \reff{ODE} we have
$$
a(x,\x,u(x),u'(x))-a_0(x,u_0(x),u_0'(x))
=-\int_x^1\left(f(y,y/\ve,u(y))-\int_0^1f(y,z,u_0(y))dz\right)dy.
$$
Hence, Lemma \ref{prep1} below yields that 
$\|a(x,\x,u(x),u'(x))-a_0(x,u_0(x),u_0'(x))\|$ is small uniformly with respect to $x$ if $\ve$ and $\|u-u_0\|_\infty$ are small.
Therefore, in this case we have local uniqueness of weak solutions $u$ to \reff{ODE} with $\ve+\|u-u_0\|_\infty\approx 0$.
\end{remark}

\begin{remark}
Our notion of weak solutions $u$ to \reff{ODE}-\reff{bc} does not require
$u \in W^{1,2}((0,1);\R^n)$, as usual,
but 
$u \in C^1([0,1];\R^n)$.
We do this in order to avoid any growth restrictions for the functions $f(x,y,u,\cdot)$. If we would suppose that the functions $f(x,y,u,\cdot)$ have linear growth
(additionally to the growth restriction for $a(x,y,u,\cdot)$, which is included in \reff{mon}), then any solution $u \in W^{1,2}((0,1);\R^n)$ to \reff{vareq} with $u(0)=0$ would be $C^1$-smooth and, hence, a weak solution to \reff{ODE}-\reff{bc} in the sense introduced above. But for that we have to use assumption \reff{mon} (cf. the proof of Lemma \ref{prep} below).

We do not believe that the assertions of Theorem \ref{main} remain true, in general, if one works with the $W^{1,2}$-notion of  weak solutions and if the maps $a(x,y,u,\cdot)$ are not monotone (even if these maps as well as the maps 
$f(x,y,u,\cdot)$
are smooth and bounded). The reason for that conjecture are the counterexamples in \cite{Ball}, which show that weak solutions $u=u_0\in W^{1,2}(0,1)$ to Dirichlet problems for ODEs of the type $a(u'(x))'=0$ with bounded and smooth, but non-monotone functions $a:\R\to \R$ are not isolated, in general, even if $a'(u'_0(x))\ge \mbox{const}>0$ for almost all $x\in(0,1)$.
\end{remark}

\begin{remark}
Theorem \ref{main} concerns only those
weak solution to \reff{hombvp} which satisfy the condition $\|u_0\|_\infty<1$. This is because in assumption \reff{mon}
the maps $a(x,y,u,\cdot)$ are supposed to be strongly monotone and Lipschitz continuous 
uniformly with respect to $x\in [0,1]$
and $y\in \R$, but to $u\in \R^n$ with $\|u\|\le 1$ only, and, hence, because the homogenized maps $a_0$ and $b_0$
are defined for $u\in \R^n$ with $\|u\|\le 1$ only, in general. Of course, this can be generalized by supposing that the maps $a(x,y,u,\cdot)$ are strongly monotone and Lipschitz continuous 
uniformly with respect to $x\in [0,1]$, $y\in \R$ and $u\in \R^n$ with $\|u\|\le R$ and with sufficiently large $R>0$.
\end{remark}

\begin{remark}
The homogenized version $a_0$ of the map $a$ depends on $a$ only (cf. \reff{bdef} and
\reff{anulldef}), but the homogenized version $f_0$ of the map $f$ depends not only on $f$, but also on $a$ (cf. \reff{fnulldef}), i.e. the homogenization of $f$ is "relative to $a$". For linear problems this effect is well-known, cf.   \cite[Remark 1.13.1]{Ben}, \cite {QXu} and \cite[formula (3.9)]{Xu}.
\end{remark}

\begin{remark}
It is easy to verify that the map $b$, which is defined in \reff{bdef},
is continuous and its first partial derivatives with respect to 
the third and fourth arguments exist and
are continuous, that the maps $b(x,\cdot,u,u')$ are 1-periodic, and that
\begin{eqnarray}
\label{bnullmon}
\big(b(x,y,u,u_1')-
b(x,y,u,u_2')\big)\cdot(u_1'-u_2') \ge 
\frac{m}{M^2}\|u_1'-u_2'\|^2,\nonumber\\
\|b(x,y,u,u_1')-b(x,y,u,u_2')
\|\le \frac{1}{m}\|u'_1-u'_2\|
\label{bnullLip}
\end{eqnarray}
for all 
$x \in [0,1]$, $y \in \R$
and $u,u_1',u_2' \in \R^n$
with $\|u\|\le 1$.
Similarly, also the map $a_0$, which is defined in \reff{anulldef},
is continuous and its first partial derivatives with respect to 
the third and fourth arguments exist and
are continuous,  and 
\begin{eqnarray*}
(a_0(x,u,u'_1)-a_0(x,u,u'_2))\cdot (u'_1-u'_1) &\ge&\frac{m^3}{M^2}\|u'_1-u'_2\|^2,\\
\|a_0(x,u,u'_1)-a_0(x,u,u'_2)\|&\le& \frac{M^2}{m}
\|u'_1-u'_2\|
%\|f_0(x,u,u'_1)-f_0(x,u,u'_2)\|&\le& \frac{M}{m}
%\|u'_1-u'_2\|
\end{eqnarray*}
for all $x \in [0,1]$, $u,u'_1,u'_2  \in \R^n$ with $\|u\|\le1$.
\end{remark}

\begin{remark}
In many applications the maps $a(x,y,u,\cdot)$ and
$f(x,y,u,\cdot)$ 
are affine, i.e. 
$$
a(x,y,u,u')=A(x,y,u)u'+\bar a(x,y,u)
\mbox{ and }
f(x,y,u,u')=F(x,y,u)u'+\bar f(x,y,u)
$$ 
with $n\times n$-matrices $A(x,y,u)$
and $F(x,y,u)$
and vectors $\bar a(x,y,u)$
and $\bar f(x,y,u)$. 
%Then 
%it follows that $b(x,y,u,u')=A(x,y,u)^{-1}(u'-\bar a(x,y,u))$,
%i.e.
%\begin{eqnarray*}
%&&b_0(x,u,u')=\int_0^1A(x,y,u)^{-1}(u'-\bar a(x,y,u))dy,\\
%&&a_0(x,u,u')=\left(\int_0^1A(x,y,u)^{-1}dy\right)^{-1}\left(u'+\int_0^1A(x,y,u)^{-1}\bar a%(x,y,u)dy\right).%
%\end{eqnarray*}
%Hence, 
%\begin{eqnarray*}
%&&f(x,y,u,b(x,y,u,a_0(x,y,u,u')))\\
%&&=\bar f(x,y,u)+F(x,y,u)A(x,y,u)^{-1}
%\left(\int_0^1A(x,z,u)^{-1}dz\right)^{-1}u'\\
%&&\;\;\;\;+F(x,y,u)A(x,y,u)^{-1}\left(\left(\int_0^1A(x,z,u)^{-1}dz\right)^{-1}
%\int_0^1A(x,z,u)^{-1}\bar a(x,z,u)dz
%-\bar a(x,y,u)\right).
%\end{eqnarray*}
%Therefore
Then also the maps $a_0(x,u,\cdot)$ and
$f_0(x,u,\cdot)$, defined by \reff{anulldef} and \reff{fnulldef}, 
are affine, i.e. 
$$
a_0(x,u,u')=A_0(x,u)u'+\bar a_0(x,u)
\mbox{ and }
f_0(x,u,u')=F_0(x,u)u'+\bar f_0(x,u)
$$ 
with
\begin{eqnarray*}
A_0(x,u)&:=&\left(\int_0^1 A(x,y,u)^{-1}dy\right)^{-1},\\
\bar a_0(x,u)&:=&\left(\int_0^1 A(x,y,u)^{-1}dy\right)^{-1}\int_0^1A(x,y,u)^{-1}\bar a(x,y,u)dy,\\
F_0(x,u)&:=&\int_0^1 F(x,y,u)A(x,y,u)^{-1}dy
\left(\int_0^1A(x,y,u)^{-1}dy\right)^{-1}
,\\
\bar f_0(x,u)&:=&\int_0^1\bar f(x,y,u)dy\\
&&+\int_0^1F(x,y,u)A(x,y,u)^{-1}dy
\left(\int_0^1A(x,y,u)^{-1}dy\right)^{-1}\int_0^1A(x,y,u)^{-1}\bar a(x,y,u)dy\\
&&-\int_0^1F(x,y,u)A(x,y,u)^{-1}\bar a(x,y,u)dy.
%-F(x,y,u)A(x,y,u)^{-1}\bar a(x,y,u)\right)dy.
\end{eqnarray*}
In particular, if $\bar a$ does not depend on $y$, then 
$\bar a_0=\bar a$
and $\bar f_0(x,u)=\int_0^1\bar f(x,y,u)dy$.
And if $F$ does not depend on $y$,  
then $F_0=F$.
%Remark that $A_0$ depends on $A$ only, but 
%$\bar a_0$ depends on $\bar a$ and $A$,
%$F_0$ depends on $F$ and $A$, and $\bar f_0$ depends on $\bar f$, $A$, $\bar a$ and $F$.
\end{remark}

\begin{remark}
\label{sufficient}
The assumption of Theorem \ref{main}, that there do not exist nontrivial weak solutions to \reff{hombvp}, is rather implicit. But there exist simple explicit sufficient contitions for it. For example, if not only the matrices $\partial_{u'}a(x,y,u_0(x),u'_0(x))$ are positive definit
(this follows from assumption \reff{mon}), but also the matrices $\partial_{u}f(x,y,u_0(x),u'_0(x))$, and if the corresponding definitness coefficients are sufficiently large in comparison with the matrix norms of $\partial_{u}a(x,y,u_0(x),u'_0(x))$ and $\partial_{u'}f(x,y,u_0(x),u'_0(x))$, then there do not exist nontrivial weak solutions to \reff{hombvp}.
\end{remark}

\begin{remark}
The assertions of Theorem \ref{main} remain true also in cases where the maps $a(x,\cdot,u,u')$ or $f(x,\cdot,u,u')$ are allowed to be discontinous, for example,
if
$$
a(x,y,u,u')=a_1(x,u,u') a_2(y)
\mbox{ and }
f(x,y,u,u')=f_1(x,u,u')f_2(y)
$$
with vector functions $a_1,f_1 \in C^1([0,1]\times \R^n\times \R^n;\R^n)$ and 1-periodic functions $a_2,f_2 \in L^\infty(\R)$, or, more general, if the maps $(x,u,u') \mapsto a(x,\cdot,u,u')$ and $(x,u,u') \mapsto f(x,\cdot,u,u')$ are continuous from $[0,1]\times \R^n\times\R^n$ into $L^\infty(\R)$
(cf. also Remark \ref{inf}). For the case of linear scalar equations see \cite[Theorems 6.1 and 6.3]{Ben} and \cite[Theorem 1.2]{YXu}.
\end{remark}

\begin{remark}
\label{Neuk}

What concerns existence and local uniqueness
for nonlinear periodic homogenization problems
(without assumption of global uniqueness) we know only the result 
\cite{Bun}  for scalar semilinear elliptic PDEs of the type
$
\mbox{\rm div}\, a(\x) \nabla u(x)=f(x)g(u(x)),
$
where the nonlinearity $g$ is supposed to have a sufficiently small local Lipschitz constant (on an appropriate bounded interval). Let us mention also \cite{Lanza1,Lanza2}, where existence and local uniqueness for a homogenization problem for the linear Poisson equation with periodic nonlinear Robin boundary conditions is shown. There the specific structure of the problem (no highly oscillating diffusion coefficients) allows to apply the classical implicit function theorem.

$L^\infty$-estimates of the homogenization error $u_\ve-u_0$ exist, to the best of our knowledge, for linear periodic homogenization problems only: For  scalar ODEs
of the type
$
\big(a(\x)u'(x)\big)'=f(x) 
$
(with a smooth 1-periodic function $a:\R \to \R$ and a smooth
function $f:[0,1] \to \R$)
in \cite[Section 1]{N}, for scalar ODEs with stratified structure of the type
$
\big(a(x,\rho(x)/\ve)u'(x)\big)'=f(x)
$
in \cite[Theorem 1.2]{YXu}.
For $L^\infty$ homogenization error estimates for scalar linear elliptic PDEs of the type
$
\mbox{\rm div}\, a(\x) \nabla u(x)=f(x)
$
see, e.g. \cite[Chapter 2.4]{Ben} and \cite{He} and for linear elliptic systems \cite[Theorem 7.5.1]{Shen}.
%(there the case of smooth coefficients is considered, and $O(\ve)$ rates are shown) and \cite{He} (there the case of two space dimensions and of non-smooth coefficients is considered, and $O(\ve)|\ln \ve|^2$ rates are shown).
\end{remark}

\begin{remark}
Consider the quasilinear elliptic PDE
$\nabla_x\cdot a(x,\x,u(x),\nabla_x u(x))=f(x)
%\mbox{\rm div}\,a(x,\x,u(x),\nabla u(x))=f(x)
% \mbox{ for } x \in \Omega \subseteq \R^d
$
with flux function $a:\Omega \times \R^d\times \R\times \R^d \to \R^d$ 
(with $\Omega \subseteq \R^d$)
such that $a(x,\cdot,u,v)$ is $\Z^d$-periodic
%$a(x,y+e_j,u,v)=a(x,y,u,v)$ for all $j=1,\ldots,d$ ($e_1:=(1,0,\ldots,0,0),\ldots,e_d:=(0,0,\ldots,0,1)$ is the standard basis in $\R^d$)
and that $a(x,y,u,\cdot)$ is strongly monotone and Lipschitz continuous uniformly with respect to $x$, $y$ and bounded $u$.
The usual formula for the homogenized flux function is (cf., e.g. \cite{Car,Fusco,Past,Wang})
\bee
\label{PDEanulldef}
a_0(x,u,v):=\int_{[0,1]^d}a(x,y,u,v+\nabla_yw(x,y,u,v))dy,
\ee
where the $\Z^d$-periodic corrector $w:\R^d\to \R$ is the solution (which depends parametrically on $x\in\Omega$, $u\in \R$ and $v\in \R^d$) of the cell problem
$\nabla_y\cdot a(x,y,u,v+\nabla_yw(y))=0$,
%$w(x,y+e_j,u,v)=w(x,y,u,v)$ 
%for $y \in \R^d$, $j=1,\ldots,d$, 
$\int_{(0,1)^d}w(y)dy=0$.
%\begin{eqnarray*}
%&&\mbox{div}_y\;a(x,y,u,v+\nabla_yw(x,y,u,v))=0 \mbox{ for } y \in \R^d,\\
%&&w(x,y+e_j,u,v)=w(x,y,u,v) \mbox{ for } j=1,\ldots,d,\;
%\int_{[0,1]^d}w(x,y,u,v)dy=0.
%\end{eqnarray*}
In space dimension one, i.e. $d=1$, this looks as follows:
\bee
\label{dual}
a_0(x,u,v):=\int_0^1a(x,y,u,v+\partial_yw(x,y,u,v))dy
\ee
and
\bee
\label{syst}
\left.
\begin{array}{l}
\displaystyle\frac{d}{dy}
a(x,y,u,v+\partial_yw(x,y,u,v))=0,\\ 
%\mbox{ for } y \in \R,\\
w(x,y+1,u,v)=w(x,y,u,v),\;
\displaystyle\int_0^1w(x,y,u,v)dy=0.
\end{array}
\right\}
\ee
From \reff{syst} it follows that $a(x,y,u,v+\partial_yw(x,y,u,v))$ is constant with respect to $y$. Therefore \reff{dual}
yields that
$
a(x,y,u,v+\partial_yw(x,y,u,v))=
a_0(x,u,v)
$ 
and, hence, 
%\label{abequ}
$
b(x,y,u,a_0(x,u,v))=v+\partial_yw(x,y,u,v),
$
i.e.
%From \reff{abequ} follows
$$
v=\int_0^1(v+\partial_yw(x,y,u,v))dy=
\int_0^1b(x,y,u,a_0(x,u,v))dy=
b_0(x,u,a_0(x,u,v)).
$$
On the other hand, the solution to \reff{syst} with $v=b_0(x,u,\bar v)$ (with arbitrary $\bar v \in \R$) is
\begin{eqnarray*}
&&w(x,y,u,b_0(x,u,\bar v))\\
&&=
\int_0^y(b(x,z,u,\bar v)-b_0(x,u,\bar v))dz-
\int_0^1\int_0^{z_1}(b(x,z_2,u,\bar v)-b_0(x,u,\bar v))dz_2dz_1.
\end{eqnarray*}
Therefore
\begin{eqnarray*}
a_0(x,u,b_0(x,u,\bar v))&=&
\int_0^1a(x,y,u,b_0(x,u,\bar v)-\partial_yw(x,y,u,b_0(x,u,\bar v)))dy\\
&=&
\int_0^1a(x,y,u,b(x,y,u,\bar v)dy=\bar v.
\end{eqnarray*}
It follows that $a_0(x,u,\cdot)=b_0(x,u,\cdot)^{-1}$.
In other words: Our definition \reff{anulldef} of the homogenized flux function $a_0$ is the same as  \reff{PDEanulldef}, considered in the  case $d=1$. In the linear case this has been shown in \cite[Remark 2.3]{Ben} and \cite[Proposition 6.16]{Ci}. 
%In \cite[Remark 5.9]{Ben} the formulas \reff{anulldef} and \reff{dual} are called dual formulas (there the linear case is considered, but with multidimensional space variable $x$).
\end{remark}

~\\ 

Our paper is organized as follows: 
In Section \ref{secabstract} we consider abstract nonlinear parameter depending equations of the type
\bee
\label{intrabstract}
\F_\ve(w)=0.
\ee
Here $\ve>0$ is the parameter. We prove a result on existence and local uniqueness of a family of solutions $w=w_\ve \approx w_0$ to \reff{intrabstract} with $\ve \approx 0$, where $w_0$ is an approximate solution to \reff{intrabstract}, i.e. an element with
$\F_\ve(w_0)\to 0$ for $\ve \to 0$, and we estimate the norm of the error $w_\ve-w_0$ by the norm of the discrepancy $\F_\ve(w_0)$. 
This type of generalized implicit function theorems has been successfully applied to singularly perturbed ODEs and PDEs in \cite{Butetc,But2022,Fiedler,NURS,
OmelchenkoRecke2015,Recke2022,
ReckeOmelchenko2008}). 
Contrary to the classical implicit function theorem it is not supposed that
the linearized operators $\F'_\ve(u)$ converge for $\ve \to 0$ in the uniform operator norm. And, indeed, in the applications to singularly perturbed problems as well as to periodic homogenization problems they do not converge for $\ve \to 0$ in the uniform operator norm (cf. Remark \ref{linop} below).
Hence, the present paper is a first step (on the ODE level) to create a common approach to  existence, local uniqueness and error estimates for singularly perturbed problems and for homogenization problems. In Part II \cite{II} we apply this approach to periodic homogenization for semilinear elliptic PDE systems.
%For another result of this type see \cite[Theorem 2.1]{Breden}.

In Section \ref{sec3} we prove Theorem \ref{main} by means of the results of Section \ref{secabstract}. For that reason we transform the boundary value problem 
\reff{ODE}-\reff{bc} into the system 
\reff{ueq}-\reff{veq}
of integral equations, and for that system 
of integral equations
we introduce an abstract setting
of the type
\reff{intrabstract}.
For that abstract setting we have to verify the key assumptions \reff{limass} and \reff{infass} of Theorem \ref{ift}, and we do this in the Subsections \ref{subsec:Feps} and \ref{subsec:infass}, respectively. 

In order to verify assumption \reff{infass} of Theorem \ref{ift} we use again approaches which are well-known in singular perturbation theory. More exactly,
proofs of $\ve$-uniform  coercivity estimates of the type \reff{infass} for various  singularly perturbed linear differential operators
in various function spaces by assuming the contrary have been used for a long time, see, e.g., \cite[Proposition 4.1]{Casteras}, \cite[Proposition 1.2]{del},
\cite[Proposition 5.1(ii)]{HSak}, \cite[Lemma 1.3]{M1}, \cite[Lemma 1]{Taniguchi}, \cite[Proposition 7.1]{Wei}.

Finally, in Section \ref{sec:other} we show that
Theorem \ref{main}
remains true also for inhomogeneous natural boundary  conditions and for two Dirichlet boundary conditions.

%, but not, in general, for inhomogeneous Neumann boundary conditions. Remark that the difficulties with inhomogeneous Neumann boundary conditions are well-known already for scalar linear problems (see, e.g. \cite[Remark 1.2.10 and Section 1.7.1]{Ben} and \cite{Suslina}). Further, we show how to prove that the assertions of Theorem \ref{main} are true also for two Dirichlet boundary conditions.

\section{An abstract result of implicit function theorem type}
\label{secabstract}
\setcounter{equation}{0}
\setcounter{theorem}{0}
In this section we formulate and prove Theorem \ref{ift} below. 
\begin{theorem}
\label{ift}
Let be given a Banach space
$W$ with norm $\|\cdot\|_W$, an open set $W_0 \subseteq W$, an element $w_0 \in W_0$ and a family of $C^1$-maps $\F_\ve:W_0 \to W$ with $\ve>0$ as family parameter. Suppose that
\bee
\label{limass}
\|\F_\ve(w_0)\|_W\to 0 \mbox{ for } \ve \to 0.
\ee
Further, suppose that there exists $\ve_0>0$ such that
\begin{eqnarray}
\label{Fredass}
&&\F'_\ve(w_0) \mbox{ is Fredholm of index zero
from $W$ into $W$ for all $\ve \in (0,\ve_0]$},\\
\label{infass}
&&\inf
\{\|\F'_\ve(w_0)w\|_W:\; \ve \in (0,\ve_0],\;
%%(\ve,w) \in (0,\ve_0]\times W, \;
\|w\|_W=1\}=:\al>0,\\
\label{supass}
&&
\sup_{\|w\|_W\le 1}\|(\F'_\ve(w_0+w_1)-\F_\ve'(w_0))w\|_W\to 0 \mbox{ for } \ve+\|w_1\|_W\to 0.
\end{eqnarray}
Then there exist $\ve_1 \in (0,\ve_0]$ and $\delta >0$ such that for all $\ve \in (0,\ve_1]$ 
there exists exactly one  $w = w_\ve \in W_0$ with $\F_\ve(w)=0$ and $\|w-w_0\|_W \le\delta$. Moreover,
\bee
\label{abest}
\|w_\ve-w_0\|_W \le\frac{2}{\al}\|\F_\ve(w_0)\|_W.
\ee
\end{theorem}
{\bf Proof }  
The assumptions \reff{Fredass} and \reff{infass}
imply, that for all $\ve \in (0,\ve_0)$ the operator $\F_\ve'(w_0)$ 
is an isomorphism from $W$ onto $W$
and
\bee
\label{invest}
\left\|\F_\ve'(w_0)^{-1}w\right\|_W \le \frac{1}{\al}\|w\|_W \mbox{ for all } w \in W.
\ee
Hence, the maps $\G_\ve:W_0 \to W$,
$$
\G_{\ve}(w) := w - \F_\ve'(w_0)^{-1}\F_\ve(w)
$$
are well-defined. Obviously, $w$ is a fixed point of $\G_\ve$ if and only if $\F_\ve(w)=0$.

For $r>0$ denote 
$ 
\mathbb{B}_r
:=\{w \in W:\; \|w-w_0\|_W \le r\}.
$
We are going to show that for sufficiently small $\ve>0$ and $r>0$ the map $\G_\ve$ is strictly contractive from the closed ball $\mathbb{B}_r$
into itself.

In order to verify the strict contractivity of $\G_\ve$ we take 
$\ve \in (0,\ve_0]$ and $v,w \in W_0$ and estimate as follows:
\begin{eqnarray*}
\lefteqn{
\|\G_{\ve}(v) - \G_{\ve}(w)\|_W = \left\|v-w-\F_\ve(w_0)^{-1}(\F_\ve(v)-\F_\ve(w))\right\|_W}\nonumber\\
&&=\left\|\F_\ve'(w_0)^{-1}\int_0^1\left(\F_\ve'(w_0)-\F'_{\ve}(s v + (1-s) w)\right)ds (v - w)\right\|_W\nonumber\\
&&\le \frac{1}{\al}\int_0^1
\|\left(\F_\ve'(w_0)-\F'_{\ve}(s v + (1-s) w)\right)(v - w)\|_Wds.
\end{eqnarray*}
Here we used \reff{invest}.
Because of assumption \reff{supass} there exist $\ve_1 \in (0,\ve_0]$ and $r_0>0$ such that $\mathbb{B}_{r_0} \subset W_0$ and
$
\|\left(\F_\ve'(w_0)-\F'_{\ve}(s v + (1-s) w)\right)(v - w)\|_W
\le\frac{\alpha}{2}\|v-w\|_W
$
for all $\ve \in (0,\ve_1]$, $s \in [0,1]$ and $v,w \in \mathbb{B}_{r_0}$.
Hence,
\bee
\label{half}
\|\G_{\ve}(v) - \G_{\ve}(w)\|_W
\le \frac{1}{2} \|v-w\|_W
\mbox{ for all } \ve \in (0,\ve_1] \mbox{ and } v,w \in \mathbb{B}_{r_0}.
\ee

Now, let us show that $\G_\ve$ maps  $\mathbb{B}_{r_0}$
into  $\mathbb{B}_{r_0}$ for all sufficiently small $\ve>0$.
Take  $\ve\in (0,\ve_1]$ and $w \in  \mathbb{B}_{r_0}$. Then \reff{invest} and \reff{half} imply
\begin{eqnarray*}
\label{in}
&&\left\|\G_{\ve}(w) - w_0\right\|_W \le \left\|\G_{\ve}(w) - \G_\ve(w_0)\right\|_W
+\left\|\G_{\ve}(w_0) - w_0\right\|_W\\
&&\le \frac{1}{2}\left\|w-w_0\right\|_W+\left\|\F_\ve'(w_0)^{-1}\F_\ve(w_0)\right\|_W
\le \frac{r_0}{2}+\frac{1}{\al}\left\|\F_\ve(w_0)\right\|_W.
%\le r_1. 
\end{eqnarray*}
But assumption \reff{limass} yields that, if $\ve_1$ is taken sufficiently small, for all  $\ve \in (0,\ve_1]$ we have 
$\|\F_\ve(w_0)\|_W \le \al r_0/2$.
Hence, for those $\ve$ we get $\left\|\G_{\ve}(w) - w_0\right\|_W \le r_0$.

Therefore, Banach's fixed point principle yields the following: For all $\ve \in (0,\ve_1]$ 
there exists exactly one  $w = w_\ve \in
\mathbb{B}_{r_0}$ with $\F_\ve(w)=0$.

Finally, let us prove \reff{abest}.  We take $\ve \in (0,\ve_1]$ and estimate as above:
$$
\|w_\ve-w_0\|_W \le \|\G_{\ve}(w_\ve) - \G_\ve(w_0)\|_W+\|\G_{\ve}(w) - w_0\|_W
\le \frac{1}{2}\|w_\ve-w_0\|_W+\frac{1}{\al}\|\F_\ve(w_0)\|_W.
$$
Hence, \reff{abest} is true.
\qed

\begin{remark}
If there exists a $C^1$-map $\tilde{\cal F}:[0,\infty)\times W \to W$ such that $\tilde{\cal F}(\ve,w)={\cal F}_\ve(w)$ for all
$(\ve,w) \in (0,\infty)\times W$, then Theorem \ref{ift} is just the classical implicit function theorem.
But in Theorem \ref{ift}
we do not suppose that $\F_\ve'(u)$ converges for $\ve \to 0$ with respect to the uniform operator norm. 
Remark that in the classical implicit function theorem one cannot omit, in general, the assumption, that $\F_\ve'(u)$ converges for $\ve \to 0$ with respect to the uniform operator norm (cf. \cite[Section 3.6]{Katz}).
%and
%\bee
%\label{comu}
%\lim_{\ve \to 0}\lim_{t \to 0}\frac{1}{t}\Big({\cal F}_\ve(w_0+tw)-{\cal F}_\ve(w_0)\Big)=\lim_{t \to 0}\lim_{\ve \to 0}\frac{1}{t}\Big({\cal F}_\ve(w_0+tw)-{\cal F}_\ve(w_0)\Big) \mbox{ for all } w \in W.
%\ee
%If such map $\tilde{\cal F}$ does not exist, then it may happen that \reff{comu} is not true, but anyway all assumptions and conclusions of Theorem \ref{ift} are satisfied. In Section \ref{sec3} below we describe an application of Theorem \ref{ift} such that \reff{comu} is not satisfied, in general.
\end{remark}

\begin{remark}
In \cite{Butetc,But2022,Fiedler,NURS,
OmelchenkoRecke2015,Recke2022,
ReckeOmelchenko2008}) slightly more general versions of Theorem \ref{ift} are used, i.e. those with $\F_\ve$ mapping one Banach space into another one, both with $\ve$-depending norms. Moreover, there the approximate solutions are allowed to be $\ve$-depending, i.e. to be a family of approximate solutions. Hence, these versions of Theorem \ref{ift} seem to be appropriate for applications to homogenization problems with approximate solutions defined by using correctors of first or higher-order (see, e.g. \cite{Allaire,Cher,Hawa}).

For another result of the type of Theorem \ref{ift} and its applications to semilinear elliptic PDE systems with numerically  determined approximate solutions  see \cite[Theorem 2.1]{Breden}.
\end{remark}

\section{Proof of Theorem \ref{main}}
\label{sec3}
\setcounter{equation}{0}
\setcounter{theorem}{0}
In this section we will prove Theorem~\ref{main} by means of Theorem \ref{ift}. Hence, all assumptions of Theorem~\ref{main} (i.e. 
\reff{perass}-\reff{mon}, existence of the weak solution $u=u_0$ to \reff{hombvp}, non-existence of weak solutions $u\not=0$ to \reff{linbvp})
will be supposed to be satisfied. 
%(without mentioning their use). 
At places, where we use the additional assumption \reff{diffass1} of Theorem \ref{main}(ii), we will mention this explicitely.

\subsection{Transformation of the boundary value problem \reff{ODE}-\reff{bc} into a system of integral equations} 
\label{subsec:transform}
In this subsection we transform the problem of weak solutions $u$ to the boundary value problem \reff{ODE}-\reff{bc} into 
the problem of solutions $(u,v)$ to
the system of integral equations \reff{ueq}-\reff{veq} below. 
%where $v_0$ is defined by
%\bee
%\label{vnulldef}
%v_0(x):=a_0(x,u_0(x),u_0'(x)) \mbox{ for } x \in [0,1].
%\ee
In what follows we use the notation
\bee
\label{Bdef}
\BB:=\{u \in C([0,1];\R^n):\; \|u\|_\infty <1\}.
\ee
%and Lemmas \ref{prep}, \ref{prep2} and \ref{prep1} below.

\begin{lemma}
\label{prep}
For all $\ve>0$ the following is true:

(i) If $u\in \BB$ is a weak solution to \reff{ODE}-\reff{bc} 
%then $u$ is $C^1$-smooth, 
and if $v \in C([0,1];\R^n)$ is defined by
\bee
\label{vdef1}
v(x)=a(x,\x,u(x),u'(x)) \mbox{ for } x \in [0,1],
\ee
then
\begin{eqnarray}
\label{ueq}
&&u(x)=\int_0^x b(y,\y,u(y),v(y))dy
\mbox{ for } x \in [0,1],\\
\label{veq}
&&v(x)=-\int_x^1f(y,\y,u(y),b(y,\y,u(y),v(y)))dy
\mbox{ for } x \in [0,1].
\end{eqnarray}

(ii) If $(u,v) \in \BB \times C([0,1];\R^n)$ is a solution to \reff{ueq}-\reff{veq}, then
$u$ is 
a weak solution to \reff{ODE}-\reff{bc}, and
%$v(x)=a(x,\x,u(x),u'(x))$ for  all $x \in [0,1)$.
and \reff{vdef1} is satified.
\end{lemma} 
{\bf Proof }
%Take $\ve>0$. 
(i) Let $u \in \BB$ be a weak solution to \reff{ODE}-\reff{bc}. Take an arbitrary test function $\vp \in C^1([0,1];\R^n)$ with $\vp(0)=0$.
Then \reff{vareq} yields that
\begin{eqnarray*}
0&=&\int_0^1\Big(a(x,\x,u(x),u'(x))\cdot\vp'(x)+
f(x,\x,u(x),u'(x))\cdot\vp(x)\Big)dx.\\
&=&\int_0^1\left(a(x,\x,u(x),u'(x))-
\int_0^xf(y,\y,u(y),u'(y))dy\right)\cdot\vp'(x)dx.
\end{eqnarray*}
Therefore $a(x,\x,u(x),u'(x))-
\int_0^xf(y,\y,u(y),u'(y))dy$ is constant with respect to $x$. In particular, 
the function $x\mapsto 
%\alpha(x):=
a(x,\x,u(x),u'(x))$
is $C^1$-smooth, and
\bee
\label{aeq}
a(x,\x,u(x),u'(x))'=f(x,\x,u(x),u'(x)).
\ee
%and the function $x\mapsto u'(x)=b(x,\x,u(x),\alpha(x))$ 
%is continuous.
If  $v \in C([0,1];\R^n)$ is defined by
\reff{vdef1}, then 
\bee
\label{uprimeeq}
u'(x)=b(x,\x,u(x),v(x)).
\ee
Inserting \reff{vdef1} and \reff{uprimeeq} into \reff{aeq} we get \reff{veq}. Moreover, the boundary condition $u(0)=0$ and \reff{uprimeeq} yield \reff{ueq}.

(ii) Let $(u,v) \in \BB\times C([0,1];\R^n)$ be a  solution to \reff{ueq}-\reff{veq}.
From \reff{ueq} follows $u(0)=0$ and \reff{uprimeeq}, and from \reff{uprimeeq}
follows \reff{vdef1}. Therefore \reff{veq} implies  $v(1)=a(1,1/\ve,u(1),u'(1))=0$ and 
$
v'(x)=a(x,\x,u(x),u'(x))'=f(x,\x,u(x),b(x,\x,u(x),v(x))).
$
If we multiply this scalarly  by
an arbitrary test function $\vp \in C^1([0,1];\R^n)$ with $\vp(0)=0$, integrate with respect to $x$ and use the boundary condition $v(1)=0$, then we get \reff{vareq}.
\qed\\

Similarly to Lemma \ref{prep} we get the following: If a function $u$ is a weak solution to the homogenized problem \reff{hombvp} and if the function
$v \in C([0,1];\R^n)$ is defined by
$
v(x):=a_0(x,u(x),u'(x)),
$
then the pair $(u,v)$ satisfies the system of integral equations
\bee
\label{transhom}
\left.
\begin{array}{l}
\displaystyle u(x)=
\int_0^x b_0(y,u(y),v(y))dy,\\
%=\int_0^x \int_0^1b(y,z,u(y),v(y))dzdy,\\
\displaystyle v(x)=
-\int_x^1                                   
f_0(y,u(y),b_0(y,u(y),v(y)))dy
%=-\int_x^1                                   \int_0^1
%f(y,z,u(y),b(y,z,u(y),v(y)))dzdy
\end{array}
\right\}\mbox{ for } x \in [0,1]
\ee
and vice versa.
Using
the definitions \reff{bdef}-\reff{fnulldef} of
$b_0$,
$a_0$ and $f_0$ and that
$
f_0(y,u,b_0(y,u,u'))=
\int_0^1f(y,z,u,
b(y,z,u,a_0(y,u,b_0(y,u,u'))))dz=
\int_0^1f(y,z,u,
b(y,z,u,u'))dz$,
we get that system \reff{transhom} is equivalent to 
\bee
\label{uvnulleq}
\left.
\begin{array}{l}
\displaystyle u(x)=
\int_0^x \int_0^1b(y,z,u(y),v(y))dzdy,\\
\displaystyle v(x)=
 -\int_x^1                                   \int_0^1
f(y,z,u(y),b(y,z,u(y),v(y)))dzdy
\end{array}
\right\}\mbox{ for } x \in [0,1].
\ee
In particular, the function $u_0$, which is by assumption of Theorem \ref{main} a weak solution to \reff{hombvp}, and the function $v_0 \in C([0,1];\R^n)$, which is defined by
\bee
\label{vnulldef}
v_0(x):=a_0(x,u_0(x),u'_0(x)),
\ee
satisfy 
\reff{uvnulleq}.
This will be used in Subsection \ref{subsec:Feps} below.

\subsection{Transformation of the linearized homogenized boundary value problem \reff{linbvp} into a system of integral equations} 
\label{subsec:translinhom}
Let us apply Lemma \ref{prep} to the linearized homogenized boundary value problem \reff{linbvp}.
Then the roles of the nonlinear maps $a(x,y,\cdot,\cdot)$ and $f(x,y,\cdot,\cdot)$ are played by the linear maps
$$
(u,u')\mapsto \partial_{u'}a_0(x,u_0(x),u_0'(x))u'+\partial_{u}a_0(x,u_0(x),u_0'(x))u
$$
and 
$$
(u,u')\mapsto \partial_{u'}f_0(x,u_0(x),u_0'(x))u'+\partial_{u}f_0(x,u_0(x),u_0'(x))u
$$
Hence, the roles of the nonlinear maps $(u,u')\mapsto b(x,y,u,u')$ and $(u,u')\mapsto f(x,y,u,b(x,y,u,u'))$ are played by the linear maps
$$
(u,u')\mapsto \partial_{u'}a_0(x,u_0(x),u_0'(x))^{-1}(u'-\partial_{u}a_0(x,u_0(x),u_0'(x))u)
$$
and 
\begin{eqnarray*}
&&(u,u')\mapsto 
\partial_{u'}f_0(x,u_0(x),u_0'(x))
\partial_{u'}a_0(x,u_0(x),u_0'(x))^{-1}(u'-\partial_{u}a_0(x,u_0(x),u_0'(x))u)\\
&&\;\;\;\;\;\;\;\;\;\;\;\;\;\;+\partial_{u}f_0(x,u_0(x),u_0'(x))u.
\end{eqnarray*}
Therefore now the role of the system \reff{ueq}-\reff{veq}
of nonlinear integral equations is played by 
the system of linear integral equations
\bee
\label{syshomlin}
\left.
\begin{array}{l}
\displaystyle u(x)=
\int_0^x \Big(A_{11}(y)u(y)+A_{12}(y)v(y)\Big)dy,\\
\displaystyle v(x)=
-\int_x^1
\Big(A_{21}(y)u(y)+A_{22}(y)v(y)\Big)dy
\end{array}
\right\}
\mbox{ for } x \in [0,1]
\ee
with
\begin{eqnarray*}
&&A_{11}(x):=-\partial_{u'}a_0(x,u_0(x),
u'_0(x))^{-1}\partial_{u}a_0(x,u_0(x),
u'_0(x))
,\\
%\partial_{u'}a_0(x,u_0(x),
%u'_0(x))^{-1}
%,\\
&&A_{12}(x):=\partial_{u'}a_0(x,u_0(x),
u'_0(x))^{-1}
,\\
%-\partial_{u'}a_0(x,u_0(x),
%u'_0(x))^{-1}\partial_{u}a_0(x,u_0(x),
%u'_0(x))
%,\\
&&A_{21}(x):=\partial_{u}f_0(x,u_0(x),u_0'(x))\\
&&\;\;\;\;\;\;\;\;\;\;\;\;\;\;\;\;-
\partial_{u'}f_0(x,u_0(x),u_0'(x))
\partial_{u'}a_0(x,u_0(x),u_0'(x))^{-1}\partial_{u}a_0(x,u_0(x),u_0'(x)),\\
%\partial_{u'}f_0(x,u_0(x),u_0'(x))
%\partial_{u'}a_0(x,u_0(x),u_0'(x))^{-1},\\
&&A_{22}(x):=\partial_{u'}f_0(x,u_0(x),u_0'(x))
\partial_{u'}a_0(x,u_0(x),u_0'(x))^{-1}.
%\partial_{u}f_0(x,u_0(x),u_0'(x))\\
%&&\;\;\;\;\;\;\;\;\;\;\;\;\;\;\;\;-
%\partial_{u'}f_0(x,u_0(x),u_0'(x))
%\partial_{u'}a_0(x,u_0(x),u_0'(x))^{-1}\partial_{u}a_0(x,u_0(x),u_0'(x)).
\end{eqnarray*}
%&&+\int_0^1F^1(x,y)A^1(x,y)^{-1}dy\left(\int_0^1 A^1(x,y)^{-1}dy\right)^{-1}
%\int_0^1A^1(x,y)^{-1}A^0(x,y)dy.
In other words, \reff{syshomlin} is the  transformation
of the linearization in $u_0$ of the homogenization of \reff{ODE}-\reff{bc}.

Moreover, if $(u,v)\in C(\overline \Omega);\R^n)^2$ is a solution to \reff{syshomlin}, then $u$ 
is a weak solution to the linearized homogenized boundary value problem \reff{linbvp}.
But by assumption of Theorem \ref{main} there do not exist weak solutions $u\not=0$ to 
\reff{linbvp}. Hence,
there do not exist nontrivial solutions $(u,v)\in C(\overline \Omega);\R^n)^2$ to \reff{syshomlin}.

Remark that \reff{vnulldef} yields that
$u_0'(x)=b_0(x,u_0(x),v_0(x))$, and from
$a_0(x,u,b_0(x,u,u'))=u'$ follows
\begin{eqnarray*}
&&\partial_{u'}a_0(x,u_0(x),u_0'(x))=
\partial_{u'}b_0(x,u_0(x),v_0(x))^{-1},\\
&&\partial_{u}a_0(x,u_0(x),u_0'(x))=
-\partial_{u'}b_0(x,u_0(x),v_0(x))^{-1}
\partial_{u}b_0(x,u_0(x),v_0(x)).
\end{eqnarray*}
Therefore
\begin{eqnarray*}
&&A_{11}(x)=\partial_{u}b_0(x,u_0(x),
v_0(x)),\;
%\partial_{u'}b_0(x,u_0(x),
%v_0(x)),\;
A_{12}(x)=
\partial_{u'}b_0(x,u_0(x),v_0(x)),\\
&&A_{21}(x)=\partial_{u}f_0(x,u_0(x),b_0(x,u_0(x),v_0(x)))\\
&&\;\;\;\;\;\;\;\;\;\;\;\;\;\;\;\;+
\partial_{u'}f_0(x,u_0(x),b_0(x,u_0(x),v_0(x)))
\partial_{u}b_0(x,u_0(x),v_0(x)),\\
%\partial_{u'}b_0(x,u_0(x),v_0(x)),\\
%\partial_{u}b_0(x,u_0(x),
%v_0(x)),\\
&&A_{22}(x)=\partial_{u'}f_0(x,u_0(x),b_0(x,u_0(x),v_0(x)))
\partial_{u'}b_0(x,u_0(x),v_0(x)).
%\partial_{u'}f_0(x,u_0(x),b_0(x,u_0(x),
%v_0(x)))
%\partial_{u'}b_0(x,u_0(x),v_0(x)),\\
%&&A_{22}(x)=\partial_{u}f_0(x,u_0(x),b_0(x,u_0(x),v_0(x)))\\
%&&\;\;\;\;\;\;\;\;\;\;\;\;\;\;\;\;+
%\partial_{u'}f_0(x,u_0(x),b_0(x,u_0(x),v_0(x)))
%\partial_{u}b_0(x,u_0(x),v_0(x)).
\end{eqnarray*}
%Hence, the linear system \reff{syshomlin} is, on the one side, the transformation of
%\reff{linbvp}, i.e. the transformation
%of the linearization in $u=u_0$ of the homogenization of the boundary value problem \reff{ODE}-\reff{bc}, and it is, 
%on the other side, the linearization in $(u,v)=(u_0,v_0)$ of \reff{transhom}, i.e. the linearization in $(u,v)=(u_0,v_0)$ of the
%transformation of the homogenization of the boundary value problem \reff{ODE}-\reff{bc}.
On the other hand, the right-hand sides of \reff{transhom} equal to the
right-hand sides of
\reff{uvnulleq}. Hence, the linearization in $(u,v)=(u_0,v_0)$ of \reff{transhom} equals to the the linearization in $(u,v)=(u_0,v_0)$ of \reff{uvnulleq}. Therefore
\begin{eqnarray*}
&&A_{11}(x)=\int_0^1\partial_{u}b(x,y,u_0(x),
v_0(x))dy
,\\
&&A_{12}(x)=\int_0^1\partial_{u'}b(x,y,u_0(x),
v_0(x))dy,\\
&&A_{21}(x)=\int_0^1\Big(
\partial_{u}f(x,y,u_0(x),
v_0(x))
\\
&&\;\;\;\;\;\;\;\;\;\;\;\;\;\;\;\;+
\partial_{u'}f(x,y,u_0(x),u_0'(x))
\partial_{u}b(x,y,u_0(x),v_0(x))\Big)dy,\\
%=\int_0^1\partial_{u'}f(x,y,u_0(x),
%v_0(x))\partial_{u'}b(x,y,u_0(x),
%v_0(x))dy,\\
&&A_{22}(x)=\int_0^1\partial_{u'}f(x,y,u_0(x),
v_0(x))\partial_{u'}b(x,y,u_0(x),
v_0(x))dy.
%\int_0^1\Big(
%\partial_{u}f(x,y,u_0(x),
%v_0(x))
%\\
%&&\;\;\;\;\;\;\;\;\;\;\;\;\;\;\;\;+
%\partial_{u'}f(x,y,u_0(x),u_0'(x))
%\partial_{u}b(x,y,u_0(x),v_0(x))\Big)dy.
\end{eqnarray*}
Therefore Lemma \ref{prep1} below implies that for all $u,v \in C(\overline \Omega;\R^n)$ we have that
\bee
\label{Adefs}
\left.
\begin{array}{l}
\displaystyle \int_0^x\Big(A_{11}(y)u(y)+A_{12}(y)v(y)\Big)dy\\
=\displaystyle\lim_{\ve \to 0}\displaystyle \int_0^x\Big(\partial_{u}b(y,y/\ve,u_0(y),
v_0(y))u(y)+\partial_{u'}b(y,y/\ve,u_0(y),
v_0(y))v(y)\Big)dy,\\
\displaystyle \int_0^x\Big(A_{21}(y)u(y)+A_{22}(y)v(y)\Big)dy\\
=\displaystyle\lim_{\ve \to 0}\displaystyle 
\int_0^x\Big(\Big(
\partial_{u}f(y,y/\ve,u_0(y),
v_0(y))\\
\;\;\;\;\;\;\;\;\;\;\;\;+
\partial_{u'}f(y,y/\ve,u_0(y),
v_0(y))\partial_ub(y,y/\ve,u_0(y),v_0(y))\Big)u(y)\\
\;\;\;\;\;\;\;\;\;\;\;\;+\partial_{u'}f(y,y/\ve,u_0(y),
v_0(y))\partial_{u'}b(y,y/\ve,u_0(y),v_0(y))v(y)\Big)
dy
%\bee
%\label{Adefs}
%\left.
%\begin{array}{l}
%\displaystyle A_{11}(x)=\lim_{\ve \to 0}\partial_{u'}b(x,\x,u_0(x),
%v_0(x)),\\
%\displaystyle
%A_{12}(x)=\lim_{\ve \to 0}\partial_{u}b(x,\x,u_0(x),
%v_0(x)),\\
%\displaystyle A_{21}(x)=\lim_{\ve \to 0}\partial_{u'}f(x,\x,u_0(x),
%v_0(x))\partial_{u}b(x,\x,u_0(x),
%v_0(x)),\\
%\displaystyle A_{22}(x)=\lim_{\ve \to 0}\Big(
%\partial_{u}f(x,\x,u_0(x),
%v_0(x))
%\\
%\;\;\;\;\;\;\;\;\;\;\;\;\;\;\;\;+
%\partial_{u'}f(x,\x,u_0(x),u_0'(x))
%\partial_{u}b(x,\x,u_0(x),v_0(x))\Big).
\end{array}
\right\}
\ee
uniformly with respect to $x \in [0,1]$.
In other words: System \reff{syshomlin}, which is the {\bf transformation of the linearization in \boldmath$u_0$\unboldmath\, of the "BVP-homogenization" of \reff{ODE}-\reff{bc}}, equals
%$(u,v)$ is a solution to the homogenization of the linearization in $(u_0,v_0)$ of \reff{ueq}-\reff{veq}, i.e. if and only if 
to  
the {\bf "integral-system-homogenization" of the linearization in \boldmath$(u_0,v_0)$\unboldmath\, of the transformation of \reff{ODE}-\reff{bc}}.
This will be used in Subsection \ref{subsec:infass} below.

\subsection{Abstract setting of \reff{ueq}-\reff{veq}}
\label{subsec:absetting}

Now we are going to apply Theorem \ref{ift} in order to solve the boundary value problem \reff{ODE}-\reff{bc} with $\ve \approx 0$ and $\|u(x)-u_0(x)\|+\|a(x,\x,u(x),u'(x))-a_0(x,u_0(x),u'_0(x))\|
\approx 0$.
We introduce the setting of Theorem \ref{ift}
as follows:
\begin{eqnarray*}
&&W:=C([0,1];\R^n)^2,\; 
\|(u,v)\|_W:=\|u\|_\infty+\|v\|_\infty,\; W_0:=\BB\times C([0,1];\R^n),\; w_0:=(u_0,v_0),\\
&&\F_\ve(u,v):=(\U_\ve(u,v),\V_\ve(u,v))
\end{eqnarray*}
with
\begin{eqnarray*}
&&[\U_\ve(u,v)](x):=u(x)-\int_0^x b(y,\y,u(y),v(y))dy\\
&&[\V_\ve(u,v)](x):=v(x)+\int_x^1f(y,\y,u(y),b(y,\y,u(y),v(y)))dy.
\end{eqnarray*}
Here $\BB$ is the open ball in $C([0,1];\R^n)$, defined in \reff{Bdef}, $u_0$ is the solution to the linearized boundary value problen \reff{linbvp}, which is given by assumption of Theorem \ref{main}, and $v_0$ is defined in \reff{vnulldef}

Because of Lemma \ref{prep}
%, \ref{prep1} and \ref{prep2} 
we have the following: If $u\in \BB$ is a weak solution to \reff{ODE}-\reff{bc},
%with $\ve \approx 0$ and $\|u-u_0\|_\infty \approx 0$, 
then  
%there exists $v \in C([0,1];\R^n)$ with $\|v-v_0\|_\infty \approx 0$ such that 
$\F_\ve(u,v)=0$ with $v$ defined by \reff{vdef1}.
%\bee
%\label{vdef}
%v(x)=a(x,\x,u(x),u'(x)).
%\ee
And if $(u,v) \in W_0$ satisfies $\F_\ve(u,v)=0$, 
%with $\ve \approx 0$ and 
%$\|u-u_0\|_\infty+
%\|v-v_0\|_\infty \approx 0$,
then $u$ 
is a weak solution to \reff{ODE}-\reff{bc}, and \reff{vdef1} is true.

Hence, if all the assumptions of Theorem \ref{ift}, i.e. \reff{limass}-\reff{supass}, are satisfied in
the setting introduced above, then
Theorem \ref{ift} yields the assertions of Theorem \ref{main}(i), in particular \reff{abest} yields \reff{est}.
If, moreover,
\bee
\label{Feps1}
\mbox{assumption \reff{diffass1} implies that }\|\F_\ve(u_0,v_0)\|_W=O(\ve) \mbox{ for } \ve \to 0
\ee
in the setting introduced above,
then the assertion \reff{abest} of Theorem \ref{ift} yields also assertion \reff{est1} of Theorem \ref{main}.

Therefore, it remains to verify the assumptions  \reff{limass}-\reff{supass}
of Theorem \ref{ift} and  the assertion \reff{Feps1}
in the setting introduced above.

\subsection{Verification of 
\reff{limass} and \reff{Feps1}}
\label{subsec:Feps}

The following lemma is the only tool from classical homogenization theory which we are going to use.
For related results see, e.g. \cite[Lemma 1.1]{N}, \cite[Proposition 2.2.2]{Shen},
\cite[Lemma 3.1]{Xu}.
In Lemma \ref{prep1} below we use the following notation for maps $g \in C([0,1]\times \R;\R^n)$:
\begin{eqnarray*}
\omega_g(\ve)&:=&\sup\{\|g(x_1,y)-g(x_2,y)\|:\;
x_1,x_2 \in [0,1],\; y \in \R,\; |x_1-x_2|\le \ve\} \mbox{ for } \ve>0,\\
\|g\|_*&:=&\sup\{\|g(x,y)\|:\;
(x,y)\in [0,1]\times\R\}.
\end{eqnarray*}
 
\begin{lemma}
\label{prep1}
Let be given 
$g \in C([0,1]\times \R;\R^n)$ such that
$g(x,y+1)=g(x,y)$ for all $x \in [0,1]$ and
$y \in \R$.
Then for all $x \in [0,1]$ and $\ve>0$ we have
\bee
\label{fsup}
\left\|
\int_{0}^{x}\left(g(y,\y)-\int_0^1g(y,z)dz\right)dy\right\|\le 2\left(\omega_g(\ve)+\ve\|g\|_*\right)
\ee
and,
if the partial derivative $\partial_xg$ exists and is continuous, 
\bee
\label{fsup1}
\left\|
\int_{0}^{x}\left(g(y,\y)-\int_0^1g(y,z)dz\right)dy\right\|\le 2\ve\left(\|g\|_*+\|\partial_xg\|_*\right).
\ee
\end{lemma}
{\bf Proof }
Define
$
h(x,y):=g(x,y)-\int_0^1g(x,z)dz.
$
Then $h(x,y+1)=h(x,y)$ and $\int_y^{y+1}h(x,z)dz=0$
%for all $x \in [0,1]$ and $y \in \R$,
and $\omega_h(\ve) \le 2\omega_g(\ve)$. 
%for all $\ve>0$. 
Therefore for $x \in [0,1]$ and $\ve>0$ we have
\begin{eqnarray*}
&&\int_{0}^{x}\left(g(y,\y)-\int_0^1g(y,z)dz\right)dy=
\int_{0}^{x}h(y,\y)dy=\ve \int_0^{x/\ve}
h(\ve y,y)dy\\
&&=\ve\left(\sum_{j=1}^{[x/\ve]}\int_{j-1}^{j}h(\ve y,y)dy+
\int_{[x/\ve]}^{x/\ve}h(\ve y,y)dy
\right)\\
&&=\ve\left(\sum_{j=1}^{[x/\ve]}\int_{j-1}^{j}\left(h(\ve y,y)dy-h(\ve j,y)\right)dy+
\int_{[x/\ve]}^{x/\ve}h(\ve y,y)dy
\right),
\end{eqnarray*}
where $[x/\ve]$ is the integer part of $x/\ve$, i.e. the largest integer which is not larger than $x/\ve$. In particular, $\ve[x/\ve]\le x \le 1$.

For $y \in [j-1,j]$ we have that $0 \le\ve(j-y)\le \ve$ and, hence, that $|h(\ve y,y)-h(\ve j,y)|\le w_h(\ve)$. Therefore
$$
\left\|\int_{0}^{x}\left(g(y,\y)-\int_0^1g(y,z)dz\right)dy\right\|
\le\ve\left([x/\ve]w_h(\ve)+\|h\|_*\right)
\le2\left(w_g(\ve)+\ve\|g\|_*\right),
$$
i.e. \reff{fsup} is proved.

If $\partial_xg$ exists and is continuous, then $g(x_1,y)-g(x_2,y)=(x_1-x_2)\int
_0^1\partial_xg(sx_1+(1-s)x_2,y)ds$, i.e.
$\omega_g(\ve) \le \ve\|\partial_xg\|_*$. Hence, in that case
\reff{fsup} implies \reff{fsup1}.
\qed\\

\begin{remark}
\label{inf}
If $g(x,y)=g_1(x)g_2(y)$ with $g_1 \in C([0,1];\R^n)$ and
$g_2 \in L^\infty(\R)$ or, more general, if
the map $x\mapsto g(x,\cdot)$ is continuous from $[0,1]$ into $L^\infty(\R)$, 
then the assertions of Lemma \ref{prep1} remain true.
\end{remark}

Because 
%of \reff{bdef} and  \reff{fnulldef} and because 
$u_0$ (the solution to the homogenized boundary value problen \reff{hombvp}, which is given by assumption of Theorem \ref{main}) and $v_0$ (defined in \reff{vnulldef}) satisfy
\reff{uvnulleq}, we have
$$
[\U_\ve(u_0,v_0)](x)=\int_0^x \left(
\int_0^1b(y,z,u_0(y),v_0(y))dz-
b(y,\y,u_0(y),v_0(y))\right)dy
$$
and 
\begin{eqnarray*}
&&[\V_\ve(u_0,v_0)](x)\\
&&=\int_x^1\Big(f(y,\y,u_0(y),b(y,\y,u_0(y),v_0(y)))-\int_0^1f(y,z,u_0(y),b(y,z,u_0(y),v_0(y)))dz  
\Big)dy.
\end{eqnarray*}
Therefore, for $\ve \to 0$ we get (by applying Lemma \ref{prep1}  with
$g(x,y)=b(x,y,u_0(x),v_0(x))$
and
$g(x,y)=f(x,y,u_0(x),b(x,y,u_0(x),v_0(x)))$,
respectively), 
that
$$
\|\U_\ve(u_0,v_0)\|_\infty+
\|\V_\ve(u_0,v_0)\|_\infty=
\left\{
\begin{array}{l}
o(1),\\
O(\ve), \mbox{ if \reff{diffass1} is satified.}
\end{array}
\right.
$$

\subsection{Verification of \reff{Fredass}-\reff{supass}} 
\label{subsec:infass}
We have
$$
[\F'_\ve(u,v)](\bar u,\bar v)=(\partial_u\U_\ve(u,v)\bar u+
\partial_v\U_\ve(u,v)\bar v,
\partial_u\V_\ve(u,v)\bar u+
\partial_v\V_\ve(u,v)\bar v)
$$
with
\begin{eqnarray*}
&&[\partial_u\U_\ve(u,v)\bar u](x)=\bar u(x)-\int_0^x \partial_ub(y,\y,u(y),v(y))\bar u(y)dy,\\
&&[\partial_v\U_\ve(u,v)\bar v](x)=-\int_0^x \partial_{u'}b(y,\y,u(y),v(y))\bar v(y)dy,\\
&&[\partial_u\V_\ve(u,v)\bar u](x)=\int_x^1\Big(\partial_uf(y,\y,u(y),b(y,\y,u(y),v(y))\\
&&\;\;\;\;+
\partial_{u'}f(y,\y,u(y),b(y,\y,u(y),v(y))\partial_ub(y,\y,u(y),v(y))\Big)\bar u(y)dy,\\
&&[\partial_v\V_\ve(u,v)\bar v](x)=\bar v(x)\\
&&\;\;\;\;+\int_x^1
\partial_{u'}f(y,\y,u(y),b(y,\y,u(y),v(y))\partial_{u'}b(y,\y,u(y),v(y))\bar v(y)dy.
\end{eqnarray*}
In order to verify assumption \reff{supass}
of Theorem \ref{ift} we calculate as follows:
\begin{eqnarray*}
&&[(\partial_u\U_\ve(u_0,v_0)
-\partial_u\U_\ve(u_1,v_1))
\bar u](x)\\
&&=\int_0^x \partial_ub(y,\y,u_1(y),v_1(y))
-\partial_ub(y,\y,u_0(y),v_0(y))
\bar u(y)dy.
\end{eqnarray*}
But $\partial_ub$ is uniformly continuous on
$\{(x,y,u,v) \in [0,1]\times \R\times \R^n\times \R^n: \; \|u\|\le 1,\; \|v\|\le R\}$ for all $R>0$. Hence,
$$
\|\partial_ub(y,\y,u_0(y),v_0(y))
-\partial_ub(y,\y,u_1(y),v_1(y))
\bar u(y)\|\to 0
\mbox{ for } \|u_0-u_1\|_\infty+
\|v_0-v_1\|_\infty\to 0
$$
uniformly with respect to $\ve>0$, $y \in [0,1]$ and $\|\bar u\|_\infty \le 1$.
Similarly one can estimate the terms
$(\partial_v\U_\ve(u_0,v_0)
-\partial_u\U_\ve(u_1,v_1))
\bar v$,
$(\partial_u\V_\ve(u_0,v_0)
-\partial_u\V_\ve(u_1,v_1))
\bar u$ and  $(\partial_v\V_\ve(u_0,v_0)
-\partial_u\V_\ve(u_1,v_1))
\bar v$.
Hence, \reff{supass} is satisfied.

Further, for any $(u,v) \in \W$ we have that $\F'_\ve(u,v)-I$ ($I$ is the identity in $W$) is a linear bounded operator from $W$ into $C^1([0,1];\R^n)^2$, where $C^1([0,1];\R^n)$ is equipped with its usual norm $\|u\|_\infty+\|u'\|_\infty$. Hence, the Arcela-Ascoli Theorem yields that for any $(u,v) \in W$ the operator $\F'_\ve(u,v)-I$ is compact from $W$ into $W$, and, therefore, the operator
$\F'_\ve(u,v)$ is Fredholm of index zero from $W$ into $W$, i.e. \reff{Fredass} is satisfied.\\

Now, let us verify \reff{infass}.

Suppose that \reff{infass} is not true, i.e. that it is not 
true that there exists $\ve_0>0$ such that $\inf\{\|\F_\ve'(w_0)w\|_W:\; \ve \in (0,\ve_0], w \in W, \|w\|_W=1\}>0$.
Then there exist sequences $\ve_1,\ve_2,\ldots>0$ and $u_1,u_2,\ldots \in C([0,1];\R^n)$ and $v_1,v_2,\ldots \in C([0,1];\R^n)$ such that
\bee
\label{contlim}
\left.
\begin{array}{r}
\displaystyle\lim_{k \to \infty}\ve_k=0,\\
\displaystyle\lim_{k \to \infty}\|\partial_u\U_{\ve_k}(u_0,v_0)u_k+
\partial_v\U_{\ve_k}(u_0,v_0)v_k\|_\infty=0,\\
\displaystyle\lim_{k \to \infty}
\|\partial_u\V_{\ve_k}(u_0,v_0)u_k+
\partial_v\V_{\ve_k}(u_0,v_0)v_k\|_\infty
=0,
\end{array}
\right\}
\ee
but
\bee
\label{conteq}
\|u_k\|_\infty+\|v_k\|_\infty
=1
\mbox{ for all } k \in \N.
\ee
Denote
\begin{eqnarray*}
%\label{undef}
&&\bar u_k(x):=\int_0^x \Big(\partial_ub(y,y/\ve_n,u_0(y),v_0(y))u_k(y)+ \partial_{u'}b(y,y/\ve_n,u_0(y),v_0(y))v_k(y)\Big)dy,\\
&&\bar v_k(x):=\int_x^1\Big(\big(\partial_uf(y,y/\ve_n,u_0(y),b(y,y/\ve_n,u_0(y),v_0(y))\nonumber\\
&&\;\;\;\;\;\;\;\;\;\;\;\;+
\partial_{u'}f(y,y/\ve_n,u_0(y),b(y,y/\ve_n,u_0(y),v_0(y))\partial_ub(y,y/\ve_n,u_0(y),v_0(y))\big)u_k(y)\nonumber\\
&&\;\;\;\;\;\;\;\;\;\;\;\;+
\partial_{u'}f(y,y/\ve_n,u_0(y),b(y,y/\ve_n,u_0(y),v_0(y))\partial_{u'}b(y,y/\ve_n,u_0(y),v_0(y))v_k(y)\Big)dy\nonumber.
\end{eqnarray*}
Then
%\bee
%\label{sup1}
$\sup_{k \in \N}\big(\|\bar u'_k\|_\infty+
\|\bar v'_k\|_\infty\big)<\infty$.
Hence, because of the Arzela-Ascoli Theorem without loss of generality we may assume that there exist $\bar u_0,\bar v_0 \in C([0,1];\R^n)$ such that
\bee
\label{unlim}
\lim_{k \to \infty}\big(\|\bar u_k-\bar u_0\|_\infty+
\|\bar v_k-\bar v_0\|_\infty\big)=0.
\ee
But we have that
\begin{eqnarray*}
&&u_k=\bar u_k+\partial_u\U_{\ve_k}(u_0,v_0)u_k+
\partial_v\U_{\ve_k}(u_0,v_0)v_k,\\
&&v_k=\bar v_k-\partial_u\V_{\ve_k}(u_0,v_0)u_k+
\partial_v\V_{\ve_k}(u_0,v_0)v_k.
\end{eqnarray*}
Hence, \reff{contlim} 
and \reff{unlim} 
yield that
\bee
\label{zero}
\lim_{k \to \infty}\big(\|u_k-\bar u_0\|_\infty+
\|v_k-\bar v_0\|_\infty\big)=0,
\ee
and \reff{conteq} implies that
\bee
\label{conteq1}
\|\bar u_0\|_\infty+\|\bar v_0\|_\infty
=1.
\ee

We are going to show that \reff{contlim},
\reff{zero} 
and \reff{conteq1} lead to a contradiction.
%imply that $\bar u_0$ is a nontrivial weak solution to the linearized boundary value problem \reff{linbvp}, what contradicts to the assumption of Theorem \ref{main} that those solutions do not exist. 
Because of \reff{Adefs}, \reff{contlim} and 
\reff{zero} 
we have for all $x\in [0,1]$ that 
\begin{eqnarray*}
&&0=\lim_{k\to\infty}\Big([\partial_u\U_{\ve_k}(u_0,v_0)u_k](x)+
[\partial_v\U_{\ve_k}(u_0,v_0)v_k](x)\Big)\\
&&=\lim_{k\to\infty}\left(u_k(x)-\int_0^x\Big(\partial_ub(y,y/\ve_k,u_0(y),v_0(y))u_k(y)+\partial_{u'}b(y,y/\ve_k,u_0(y),v_0(y))v_k(y)\Big)dy\right)\\
&&=\bar u_0(x)-\lim_{k\to\infty}\left(\int_0^x\Big(\partial_ub(y,y/\ve_k,u_0(y),v_0(y))\bar u_0(y)+\partial_{u'}b(y,y/\ve_k,u_0(y),v_0(y))\bar v_0(y)\Big)dy\right)\\
&&=\bar u_0(x)-\int_0^x\Big(A_{11}(y)\bar u_0(y)
+A_{12}(y)\bar v_0(y)\Big)dy.
\end{eqnarray*}
Hence, the pair $(\bar u_0,\bar v_0)$ satisfies the first equation of system \reff{syshomlin}.

Similarly one shows that the pair $(\bar u_0,\bar v_0)$ satisfies also the second equation of system \reff{syshomlin}. Indeed, as above we have
\begin{eqnarray*}
0&=&\lim_{k\to\infty}\Big([\partial_u\V_{\ve_k}(u_0,v_0)u_k](x)+
[\partial_v\V_{\ve_k}(u_0,v_0)v_k](x)\Big)\\
&=&\lim_{k\to\infty}\left(v_k(x)+
\int_x^1\Big(\big(\partial_uf(y,y/\ve_k,u_0(y),b(y,y/\ve_k,u_0(y),v_0(y))\right.\nonumber\\
&&\;\;\;\;\;\;\;\;\;\;\;\;+
\partial_{u'}f(y,y/\ve_k,u_0(y),b(y,y/\ve_k,u_0(y),v_0(y))\partial_ub(y,y/\ve_k,u_0(y),v_0(y))\big)u_k(y)\nonumber\\
&&\;\;\;\;\;\;\;\;\;\;\;\;+
\partial_{u'}f(y,y/\ve_k,u_0(y),b(y,y/\ve_k,u_0(y),v_0(y))\partial_{u'}b(y,y/\ve_k,u_0(y),v_0(y))v_k(y)\Big)dy\Big)\\
&=&\bar v_0(x)+\lim_{k\to\infty}
\int_x^1\Big(\big(\partial_uf(y,y/\ve_k,u_0(y),b(y,y/\ve_k,u_0(y),v_0(y))\nonumber\\
&&\;\;\;\;\;\;\;\;\;\;\;\;+
\partial_{u'}f(y,y/\ve_k,u_0(y),b(y,y/\ve_k,u_0(y),v_0(y))\partial_ub(y,y/\ve_k,u_0(y),v_0(y))\big)\bar u_0(y)
\nonumber\\
&&\;\;\;\;\;\;\;\;\;\;\;\;+
\partial_{u'}f(y,y/\ve_k,u_0(y),b(y,y/\ve_k,u_0(y),v_0(y))\partial_{u'}b(y,y/\ve_k,u_0(y),v_0(y))\bar v_0(y)\Big)dy\\
&=&\bar v_0(x)+\int_x^1\Big(A_{21}(y)\bar u_0(y)
+A_{22}(y)\bar v_0(y)\Big)dy.
\end{eqnarray*}

Now we use that the system \reff{syshomlin} does not have nontrivial solutions (cf. Subsection \ref{subsec:translinhom}), i.e. $\bar u_0=\bar v_0=0$. But this is a contradiction to \reff{conteq1}.

\begin{remark}
\label{linop}
It is easy to verify that the linear operators $\F'_\ve(u,v)$ do not converge for $\ve \to 0$ in the uniform operator norm in ${\cal L}(W)$, in general. For example,
$\int_0^x\partial_ub(y,\y,u(y),v(y))\bar u(y)dy$ 
(with fixed $u,v \in C([0,1];\R^n)$)
does not converge for $\ve \to 0$ uniformly with respect to $x \in [0,1]$ and 
$\bar u\in C([0,1];\R^n)$ with $\|\bar u\|_\infty \le 1$, in general. 
\end{remark}

\section{Other boundary conditions}
\label{sec:other}
\setcounter{equation}{0}
\setcounter{theorem}{0}

\subsection{Inhomogeneous boundary conditions}
If the homogeneous natural boundary condition in \reff{bc} is replaced by a corresponding inhomogeneous one, i.e.
\bee
\label{bc1}
u(0)=0, \;a(1,1/\ve,u(1),u'(1))=u^1,
\ee
then the weak formulation of the boundary value problem \reff{ODE},\reff{bc1} is as follows: Find $u \in W^{1,2}((0,1);\R^n)$ such that $u(0)=0$ and
$$
\int_0^1\Big(a(x,\x,u(x),u'(x))\cdot\vp'(x)
+f(x,\x,u(x),u'(x))\cdot\vp(x)\Big)dx=
u^1\cdot\vp(1)
$$
for all $\vp \in C^1([0,1];\R^n)$
with $\vp(0)=0$.
The system \reff{ueq}-\reff{veq} of integral equations has to be replaced by
\begin{eqnarray*}
&&u(x)=\int_0^x 
b(y,\y,u(y),v(y))dy,\\
&&v(x)=u^1-\int_x^1f(y,\y,u(y),b(y,\y,u(y),v(y)))dy,
\end{eqnarray*}
and the results of Theorem \ref{main} remain unchanged.

\subsection{Two Dirichlet boundary conditions}
If we consider \reff{ODE} with two homogeneous Dirichlet boundary conditions, i.e.
$
u(0)=u(1)=0,
$
then the weak formulation is as follows: Find $u \in W^{1,2}((0,1);\R^n)$ such that $u(0)=u(1)=0$ and
\bee
\label{vareq2}
\left.
\begin{array}{r}
\displaystyle\int_0^1\Big(a(x,\x,u(x),u'(x))\cdot\vp'(x)
+f(x,\x,u(x),u'(x))\cdot\vp(x)\Big)dx=0\\
\mbox{ for all } \vp \in C^1([0,1];\R^n)
\mbox{ with } \vp(0)=\vp(1)=0,
\end{array}
\right\}
\ee
and the system \reff{ueq}-\reff{veq} has to be changed to the system
\begin{eqnarray*}
u(x)&=&\int_0^x 
b(y,\y,u(y),v(y))dy,\\
v(x)&=&w-\int_x^1f(y,\y,u(y),b(y,\y,u(y),v(y)))dy,\\
0&=&\int_0^1b(x,\x,u(x),v(x))dx
\end{eqnarray*}
with unknowns $(u,v,w) \in C([0,1];\R^n)\times C([0,1];\R^n) \times \R^n$.
This system of integral equations can be treated by Theorem \ref{ift} in the following setting:
$$
W:=C([0,1];\R^n)^2\times \R^n ,\; 
\|(u,v,w)\|_W:=\|u\|_\infty+\|v\|_\infty+\|w\|,\; W_0:=\BB\times C([0,1];\R^n)\times \R^n.
$$
The approximate solution $w_0$ of Theorem \ref{ift} is the triple 
$
(u_0,v_0,w_0) \in C([0,1];\R^n)^2\times \R^n,
$
 where $u_0$ is the given weak solution of the homogenized problem, i.e. 
$u_0 \in W^{1,2}((0,1);\R^n)$ with $\|u_0\|_\infty<1$ and $u_0(0)=u_0(1)=0$ and
$$
\begin{array}{r}
\displaystyle\int_0^1\Big(a_0(x,u_0(x),u_0'(x))\cdot\vp'(x)
+f_0(x,u_0(x),u_0'(x))\cdot\vp(x)\Big)dx=0\\
\mbox{ for all } \vp \in C^1([0,1];\R^n)
\mbox{ with } \vp(0)=\vp(1)=0,
\end{array}
$$
and
$
v_0(x):=a_0(x,u_0(x),u_0'(x))$, $w_0:=v_0(1)$.
The homogenized vector functions $a_0$ and $f_0$ are defined as in \reff{anulldef} and \reff{fnulldef}.
Finally, the maps $\F_\ve:W\to W$ are
defined by
$$
\F_\ve(u,v,w):=(\U_\ve(u,v),\V_\ve(u,v,w),\W_\ve(u,v))
$$
with
\begin{eqnarray*}
&&[\U_\ve(u,v)](x):=u(x)-\int_0^x b(y,\y,u(y),v(y))dy,\\
&&[\V_\ve(u,v,w)](x):=v(x)-w+\int_x^1f(y,\y,u(y),b(y,\y,u(y),v(y)))dy,\\
&&\W_\ve(u,v):=\int_0^1b(x,\x,u(x),v(x))dx.
\end{eqnarray*}
As in Subsection \ref{subsec:Feps} it follows that
\begin{eqnarray*}
&&[\U_\ve(u_0,v_0)](x)=\int_0^x \Big(
b_0(y,u_0(y),v_0(y))-
b(y,\y,u_0(y),v_0(y))\Big)dy,\\
&&[\V_\ve(u_0,v_0,w_0)](x)\\
&&=\int_x^1\Big(f(y,\y,u(y),b(y,\y,u(y),v(y)))
-f_0(y,u_0(y),b_0(y,u_0(y),v_0(y)))\Big)
dy,
\end{eqnarray*}
and, hence,  for $\ve \to 0$ follows that
$$
\|\U_\ve(u_0,v_0)\|_\infty+
\|\V_\ve(u_0,v_0,w_0)\|_\infty=
\left\{
\begin{array}{l}
o(1),\\
O(\ve), \mbox{ if \reff{diffass1} is satified.}
\end{array}
\right.
$$
Further, we have
\begin{eqnarray*}
\W_\ve(u_0,v_0)&=&\int_0^1b(x,\x,u_0(x),v_0(x))dx\\
&=&
\int_0^1(b(x,\x,u_0(x),v_0(x))
-b_0(x,u_0(x),v_0(x))
dx
\end{eqnarray*}
and, hence, for $\ve \to 0$ follows
$$
\|\W_\ve(u_0,v_0)\|_\infty=
\left\{
\begin{array}{l}
o(1),\\
O(\ve), \mbox{ if \reff{diffass1} is satified.}
\end{array}
\right.
$$
Therefore, the assumptions 
\reff{limass} and \reff{Feps1}
are satisfied in the setting introduced above.

Finally, let us verify the assumption 
\reff{infass} 
of Theorem \ref{ift} in the setting introduced above.
Suppose that \reff{infass} is not true.
Then there exist sequences $\ve_1,\ve_2,\ldots>0$ and $u_1,u_2,\ldots \in C([0,1];\R^n)$ and $v_1,v_2,\ldots \in C([0,1];\R^n)$ and $w_1,w_2,\ldots \in \R^n$ such that $\ve_n \to 0$ for $n \to \infty$ and
\begin{eqnarray}
\label{ueq1}
\lim_{k \to \infty}\|\partial_u\U_{\ve_k}(u_0,v_0)u_k+
\partial_v\U_{\ve_k}(u_0,v_0)v_k\|_\infty&=&0,\\
\label{veq1}
\lim_{k \to \infty}
\|\partial_u\V_{\ve_k}(u_0,v_0)u_k+
\partial_v\V_{\ve_k}(u_0,v_0)v_k+
\partial_w\V_{\ve_k}(u_0,v_0)w_k
\|_\infty
&=&0,\\
\lim_{k \to \infty}\|\partial_u\W_{\ve_k}(u_0,v_0)u_k+
\partial_v\W_{\ve_k}(u_0,v_0)v_k\|&=&0,\nonumber
\end{eqnarray}
but
\bee
\label{one1}
\|u_k\|_\infty+\|v_k\|_\infty+\|w_k\|
=1
\mbox{ for all } k \in \N.
\ee
As in Subsection 
\ref{subsec:infass} one can show that without loss of generality we can assume that there exist $\bar u_0, \bar v_0 \in C([0,1];\R^n)$ and $\bar w_0 \in \R^n$ such that
$
\|u_k-\bar u_0\|_\infty+
\|v_k-\bar v_0\|_\infty+
\|w_k-\bar w_0\|\to 0$ for $k \to \infty$
and that $\bar u_0$ is a solution to the linearization (in $u_0$) of \reff{vareq2}, i.e. that $\bar u_0=0$.
Hence, from \reff{ueq1} follows that for all $x \in [0,1]$ we have
\begin{eqnarray*}
0&=&\lim_{k \to \infty}[\partial_v\U_\ve(u_0,v_0)v_k](x)=
\lim_{k \to \infty}\int_0^x\partial_{u'}b(y,\y,u_0(y),v_0(y))v_k(y)dy\\
&=&
\int_0^x\partial_{u'}b_0(y,u_0(y),v_0(y))\bar v_0(y)dy.
\end{eqnarray*}
It follows that $\partial_{u'}b_0(x,u_0(x),v_0(x))\bar v_0(x)=0$ 
for all $x \in [0,1]$ and, hence, \reff{bdef} and 
%$$
%\int_0^1 \partial_{u'}b_0(x,u_0(x),v_0(x))\bar v_0(x)\cdot\
%\bar v_0(x) dx=0,
%$$
\reff{bnullmon} yield that $\bar v_0=0$.
Therefore, \reff{veq1} implies that
$$
0=\lim_{k \to \infty}\|\partial_w\V(u_0,v_0,w_0)\bar w_k\|_\infty=\lim_{k \to \infty}\|w_k\|=\|\bar w_0\|,
$$
and we get a contradiction to \reff{one1}.

\end{document}